\documentclass[11pt]{article}

\usepackage{natbib}
\usepackage{a4wide}
\usepackage{algpseudocode,algorithm,algorithmicx}
\usepackage{xcolor}
\usepackage{hyperref}

\usepackage[utf8]{inputenc}
\usepackage{longtable}
\usepackage{float}
\usepackage{amssymb}
\usepackage{amsmath}
\usepackage{adjustbox}
\usepackage{soul}
\usepackage{nomencl}
\makenomenclature
\usepackage{tablefootnote}

\usepackage{etoolbox}
\renewcommand\nomgroup[1]{%
  \item[\bfseries
  \ifstrequal{#1}{D}{Decisions variables}{%
  \ifstrequal{#1}{H}{Auxiliary variables}{%
  \ifstrequal{#1}{G}{Generator Parameters}{%
  \ifstrequal{#1}{P}{Model Generator Parameters}{%
  \ifstrequal{#1}{W}{System Parameters}{%
  \ifstrequal{#1}{Z}{Storage Parameters}{%
    \ifstrequal{#1}{F}{Functions}{%
  \ifstrequal{#1}{S}{Sets}{}}}}}}}}%
]}

\newcommand{\sUnit}{s}
\newcommand{\sUnitSet}{S}
\newcommand{\productionCharge}{pc}
\newcommand{\productionDischarge}{pd}
\newcommand{\totalEnergy}{pe}
\newcommand{\effCharge}{\eta^{c}}
\newcommand{\effDischarge}{\eta^{d}}
\newcommand{\maxCharge}{\overline{PC}}
\newcommand{\maxDischarge}{\overline{PD}}
\newcommand{\maxEnergy}{\overline{PE}}
\newcommand{\minEnergy}{\underline{PE}}

\newcommand*\rot{\rotatebox{90}}
\newcommand{\unit}{g}
\newcommand{\unitSet}{G}
\newcommand{\node}{n}
\newcommand{\nodeSet}{N}
\newcommand{\nodeInjection}{inj}
\newcommand{\timeUnit}{t}
\newcommand{\timeSet}{T}
\newcommand{\flow}{f}

\newcommand{\production}{p}

\newcommand{\productionMax}{\overline{P}}
\newcommand{\productionMin}{\underline{P}}
\newcommand{\rampUp}{RU}
\newcommand{\rampDown}{RD}

\newcommand{\demand}{D}

\newcommand{\unitCommit}{u}
\newcommand{\decisionStart}{v}

\newcommand{\decisionStop}{w}

\newcommand{\minUpTime}{UT}
\newcommand{\minDownTime}{DT}

\newcommand{\No}{{-}}
\newcommand{\Yes}{{+}}

\title{New efficient ADMM algorithm for the Unit Commitment Problem}
\author{Rogier Hans Wuijts\footnote{Utrecht University, R.H.Wuijts@uu.nl}, Marjan van den Akker\footnote{Utrecht University, J.M.vandenAkker@uu.nl}, Machteld van den Broek\footnote{University of Groningen, M.A.van.den.Broek@rug.nl}}

\begin{document}
\maketitle

\begin{abstract}
The unit commitment problem (UC) is an optimization problem concerning the operation of electrical generators. Many algorithms have been proposed for the UC and in recent years a more decentralized approach, by solving the UC with alternating direction method of multipliers (ADMM), has been investigated. For convex problems ADMM is guaranteed to find an optimal solution. However, because UC is non-convex additional steps need to be taken in order to ensure convergence to a feasible solution of high quality. Therefore, solving UC by a MIL(Q)P formulation and running an off-the-shelf solver like Gurobi until now seems to be the most efficient approach to obtain high quality solutions.

In this paper, we introduce a new and efficient way to solve the UC with ADMM to near optimality. We relax the supply-demand balance constraint and deal with the non-convexity by iteratively increasing a penalty coefficient until we eventually force convergence and feasibility. At each iteration the subproblems are solved by our efficient algorithm for the single UC subproblem developed in earlier work and our new ADMM algorithm for the transmission subproblems.

Computational experiments on benchmark instances demonstrated that our algorithm produces high-quality solutions. The computation time seems to grow practically linear with the length of the time horizon. For the case with quadratic generation cost our algorithm is significantly faster than solving the problem by a state-of-the-art MIL(Q)P formulation. For the case of linear generation cost, it outperforms the MILP approach for longer time horizons.
\end{abstract}

\thispagestyle{empty}
\nomenclature[P]{$\overline{P}_g$,$\underline{P}_g$}{maximum and minimum power output of generator $g$}
\nomenclature[D]{$u_{gt}$,$v_{gt}$,$w_{gt}$}{Binary commitment, start and stop variable of generator $g$ at time $t$.}
\nomenclature[S]{$G$,$R$,$T$,$S$,$L$}{Set of the generators, renewable generators, time steps, storage units and transmission lines}
\nomenclature[D]{$p_{gt}$}{Production of generator $g$ at time $t$}
\nomenclature[D]{$inj_{nt}$}{The amount of power node $n$ draws from the transmission system at time step $t$.}
\nomenclature[D]{$f_{lt}$}{The amount of power node that flows on the transmission line $l$. if $l = (n \rightarrow n')$ the power that flows from node $n$ to node $n'$}
\nomenclature[P]{$\overline{f_l}$,$\underline{f_l}$}{Maximum flow and maximum negative flow on transmission line $l$}
\nomenclature[P]{$\rampUp_g$,$\rampDown_g$}{Ramp up and down limit of generator $g$}
\nomenclature[P]{$SU_g$,$SD_g$}{Start-up and shutdown limit of generator $g$}
\nomenclature[P]{$UT_g$,$DT_g$}{The minimum uptime and downtime of generator $g$}
\nomenclature[P]{$\minEnergy_s$,$\maxEnergy_s$}{The minimum and maximum energy of storage unit $s$}
\nomenclature[P]{$\maxCharge_s$,$\maxDischarge_s$}{Charging and discharging limit of storage unit $s$}
\nomenclature[P]{$\effCharge_s$,$\effDischarge_s$}{Charge and discharge efficiency of storage unit $s$}
\nomenclature[D]{$\productionCharge_{\sUnit \timeUnit}$,$\productionDischarge_{\sUnit \timeUnit}$}{Charging and discharging of storage unit $s$ at time $t$}
\nomenclature[D]{ $p_s$}{Sum of charging and discharging of storage unit $s$ at time $t$}
\nomenclature[D]{$p_{rt}$}{Production of renewable energy source $r$ at time $t$}
\nomenclature[P]{$AT_{rt}$}{Availability factor of renewable source $r$ at time $t$}

\printnomenclature

\section{Introduction}
The unit commitment problem (UC) is a family of  NP-hard  \cite{van2018large,melhorn2016validating} optimization problems for the operational planning of electrical generators. The goal is to find a minimum cost schedule of generators that satisfies the electricity demand at each time step. UC is a widely studied topic in Operations Research, but its significance extends beyond this field, as it plays a central role in detailed power system modeling, which is crucial in the ongoing global shift towards sustainable energy \cite{welsch2014incorporating,abujarad2017recent}. 

The UC has been extensively researched in academia, with various optimization techniques applied to solve it. This can be seen in the numerous papers cited in various surveys on the topic \cite{van2018large,saravanan2013solution,padhy2004unit} spanning multiple algorithm paradigms, from exact methods to heuristic evolutionary algorithms. However, since UC has a mixed combinatorial and continuous statespace, mathematical programming techniques have been most successful in finding near optimal solutions \cite{van2018large}. Until the early 2000’s, Lagrangian Relaxation was the method of choice for solving UC \cite{borghetti2003lagrangian}. However, this has changed. In the past 15 years, Mixed Integer Linear Programming (MILP) problem formulations have been a very active area of research \cite{knueven2020novel} and popular tool for solving the UC \cite{atakan2017state,taktak2017overview,morales2013tight}.  Additionally, it has been stated that in the industry, UC is commonly solved using MILP and that the shift from Lagrangian Relaxation to MILP was considered beneficial \cite{o2011recent,carlson2012miso,o2007s}.

The recent  resurgence of interest in alternating direction method of multipliers (ADMM) \cite{yang2022survey} have also resulted in an increased application of augmented Lagrangian methods, such as ADMM, to the UC  \cite{feizollahi2015large,kraning2014dynamic,ramanan2017asynchronous,zhang2021computational,xavier2020decomposable}.

ADMM is an augmented Lagrangian method that solves convex optimization problems. In ADMM, decision variables are split into two or more subsets, and the coupling constraints, i.e constraints involving variables from different subsets, are relaxed by moving them to the objective such that violation of these constraints results in a penalty cost. Next, at consecutive steps subproblems are solved over a subset of the variables. These subproblems are much easier to solve than the original problem and this can be done in a distributed and parallel way. 

The UC problem can be decomposed in multiple ways. For example, ADMM can be applied to the UC by relaxing the demand coupling constraint at each node, which was done by \cite{kraning2014dynamic} and is also commonly used in existing literature on Lagrangian relaxation for UC \cite{borghetti2001lagrangian,fan2002new,frangioni2008solving}. The demand coupling constraint at a node ensures that the total production of electricity is equal to the consumption of electricity. By relaxing this constraint, the problem decomposes into multiple subproblems for the generators, Renewable Energy Sources (RES), transmission, and storage units. These subproblems can then be solved by efficient algorithms  \cite{wuijts2021improved,frangioni2006solving,guan2003optimization,FAN2002215}. The Lagrange multipliers at each node corresponds to  a nodal electricity price for each time step. 

When the problem is decomposed at the transmission level \cite{feizollahi2015large,ramanan2017asynchronous,xavier2020decomposable,zhang2021computational}, the nodes of the power system are divided into regions, which are then decoupled by copying transmission variables at their boundaries and relaxing the equality between these copied variables. If you decompose at the transmission level, the subproblems at each node are (smaller) UC problems which contain a subset of the total generators, RES and storage units. These UC problems can still be hard to solve and are often solved with MILP.

There are several other possible decompositions. For example, duplicated dispatch variables can be used in different subproblems of the UC problem \cite{batut1992daily} or temporal constraints such as ramping limits or minimum up- and downtime can be relaxed \cite{safdarian2019temporal}.

 Only for convex optimization problems, ADMM is guaranteed to find an optimal solution \cite{boyd2011distributed}. However, since UC is non-convex due to the binary commitment variables, researchers using ADMM to solve UC must be aware that it may not converge to the optimal solution, if at all. This is a general drawback of decomposition algorithms based on Lagrangian relaxation. 
 
 Some authors \cite{feizollahi2015large,xavier2020decomposable} solved the issue of nonconvexity  by embedding their ADMM algorithm into a larger heuristic that after some iterations fixes the binary commitment variables to make the problem convex. The ADMM algorithm can then converge towards the local optimal solution given these fixed binary variables. This process is called “release-and-fix”. It is not guaranteed  that these fixed commitment variables are optimal or even result in feasible solutions. 
 
 Another approach is to repair the final solution provided by Lagrangian relaxation in an ad hoc way by fixing the binary variables with some heuristic \cite{feltenmark1997optimization,yan1993scheduling,borghetti2003lagrangian,frangioni2008solving}. A convex economic dispatch problem is then solved with these fixed binary variables. For example, the heuristics proposed by  \cite{feltenmark1997optimization} sequentially and greedily processes the timesteps, fixing binary variables while adhering to minimum up times. Which generators are being turned on or off is based on either a priority list or randomly sampled from a probability distribution calculated using information
 from their bundle methods. They introduce four variants, one of which is deterministic and the others use sampling. When sampling is used, the heuristics are repeated, and the best solution is retained. The authors found that three out of the four proposed heuristics frequently fail to produce feasible solution and the fourth one requires a lot of sampling steps (repeated 200 \cite{feltenmark2000dual} or 1000 \cite{feltenmark1997optimization} times). 
 
Note that in the above heuristics ramping limits are ignored. \cite{frangioni2008solving} and \cite{borghetti2003lagrangian} both use a slight modification of one of the four heuristics proposed \cite{feltenmark1997optimization}. \cite{frangioni2008solving} studies the case of UC with ramping limits and only finds feasible solutions in roughly one out of ten dual iterations. Moreover, the authors state that incorporating ramping limits into their heuristic is non-trivial. 
 
Other authors try to force feasibility by adding an additional penalty function \cite{dubost2005primal}. However, solutions produced are still infeasible and far away from optimality. Finally, others have chosen to just ignore the issue  \cite{batut1992daily} and some authors chose to convexify the generator polytope, i.e. they remove the non-convexity by taking the convex-hull of the generator state space  \cite{kraning2014dynamic}, again leading to infeasible solutions. \\

{\bf Our contribution:} in this paper we introduce a new and efficient way to solve UC with ADMM. We relax the demand coupling constraint and improve existing approaches in two ways. Firstly, we apply a different method to deal with the non-convexity by using an increasing penalty coefficient. This means that we do not need heuristic approaches to transfer the infeasible solution into a feasible one. Our algorithm is the first Langrangian algorithm to force convergence and feasibility for the UC including ramping limits without requiring an ad-hoc repair heuristic. Secondly, we solve the subproblems very efficiently by applying our 1UC algorithm \cite{wuijts2021improved}, which strongly outperforms earlier algorithms for 1UC especially on longer time horizons, and also our newly created algorithm to solve transmission sub problems. 
We performed computational experiments on all known benchmark instances. Although our ADMM algorithm needs many iterations to converge, each iteration has a low computation time due to our efficient 1UC algorithm. 

Consequently, we have introduced a fast algorithm that produces high quality solutions respecting the ramping limits which experimentally outperforms solving UC with MIP formulations, the current state-of-the-art for the UC. Moreover, our method scales much better with the length of the planning horizon. In the experiments, we always found feasible solutions of a very good quality for a large set of UC instances. \\

The remainder of this paper is organized as follows. In the following section, we present the problem formulation of UC. In Section 3, we describe how ADMM can be applied to a small power system. In Section 4, we describe how we applied ADMM to a full UC formulation, how we solve the individual subproblems and how we dealt with non-convexity. Our experiments and results are presented in Section 5, and we end the paper with a discussion in Section 6 and conclusions in Section 7.

\section{UC problem description}\label{UC_problem_description}
In this section, we define the UC that we study in this paper. Our UC description is based on a Mixed Integer Quadratic Program (MIQP)  formulation with detailed thermal generators, RES, storage, and transmission lines.
The decision variables $p_{gt}$, $p_{rt}$, $p_{st}$ are the generation of thermal generator unit $g$, RES unit $s$ and storage unit $s$ at time step $t$. Other variables are created to specify the feasible state space of these variables. We use the well-known 3-bin formulation \cite{knueven2020novel}. This means that we use binary variables $u_{gt}$, $v_{gt}$, and $w_{gt}$ to signal that generator $g$ respectively is on, switched on, and switched off in time step $t$ . The formulation is as follows:
\begin{subequations}
\begin{align}
& \min \sum_{t \in T}\sum_{g \in G} a_g u_{gt} + b_g p_{gt} + c_g p_{gt}^2 + v_{gt} cost_{start} \label{eq:chapter3:objective}\\
& s.t. \nonumber \\
& \production_{\unit \timeUnit} \geq \unitCommit_{\unit \timeUnit} \productionMin_\unit , \unit \in \unitSet, \timeUnit \in \timeSet \label{eq:chapter3:pmin}\\
&\production_{\unit \timeUnit} \leq \productionMax_g \unitCommit_{\unit \timeUnit} , \unit \in \unitSet, \timeUnit \in \timeSet \label{eq:chapter3:pmax}\\
&\sum_{i = \timeUnit - \minUpTime_g +1}^{\timeUnit} \decisionStart_{gi} \leq \unitCommit_{\unit \timeUnit} , \timeUnit \in \timeSet, \unit \in \unitSet \label{eq:chapter3:MDT1} \\ 
&\sum_{i = \timeUnit - \minDownTime_g +1}^{\timeUnit} \decisionStop_{gi} \leq 1- \unitCommit_{\unit \timeUnit} , \timeUnit \in \timeSet, \unit \in \unitSet \label{eq:chapter3:MDT2}\\ 
&\production_{\unit \timeUnit} - \production_{\unit \timeUnit-1} \leq (SU_g - \rampUp_g) \decisionStart_{g\timeUnit} + \rampUp_g \unitCommit_{g \timeUnit} , \timeUnit \geq 2, \unit \in \unitSet \label{eq:chapter3:tightRamp1} \\ 
&\production_{\unit \timeUnit-1} - \production_{\unit \timeUnit} \leq (SD_g - \rampDown_g) \decisionStop_{g\timeUnit} + \rampDown_g \unitCommit_{g\timeUnit-1} , \timeUnit \geq 2, \unit \in \unitSet \label{eq:chapter3:tightRamp2} \\
& p_{rt} \leq AF_{rt} \productionMax_{rt} , r \in R, t \in T \label{eq:chapter3:RES}\\
& 0 \leq \productionCharge_{\sUnit \timeUnit} \leq \maxCharge_\sUnit , \timeUnit \in \timeSet, s \in S\label{eq:chapter3:ChargeLimit} \\
& 0 \leq \productionDischarge_{\sUnit \timeUnit} \leq \maxDischarge_\sUnit , \timeUnit \in \timeSet, s \in S \label{eq:chapter3:DischargeLimit}\\ 
& p_{st} = \productionDischarge_{\sUnit \timeUnit} - \productionCharge_{\sUnit \timeUnit} , \timeUnit \in \timeSet, s \in S\label{eq:chapter3:StorageProduction} \\
& \minEnergy_\sUnit \leq \totalEnergy_{\sUnit \timeUnit} \leq \maxEnergy_\sUnit , \timeUnit \in \timeSet, s \in S \label{eq:chapter3:EnergyLimit}\\
& \totalEnergy_{\sUnit \timeUnit} = \totalEnergy_{\sUnit \timeUnit - 1} + \productionCharge_{\sUnit \timeUnit} * \effCharge_{\sUnit \timeUnit} - \frac{\productionDischarge_{\sUnit \timeUnit}}{\effDischarge_{\sUnit \timeUnit}} , \timeUnit \in \timeSet, s \in S \label{eq:chapter3:NetStorageProduction}\\
& \nodeInjection_{\node t } = \sum_{l=(n' \rightarrow n), \node' \in N} \flow_{lt} , \timeUnit \in \timeSet, \node \in \nodeSet \label{eq:chapter3:trans1}\\
&\underline{\flow_{l}} \leq \flow_{lt} \leq \overline{\flow_{l}} , l \in L, \timeUnit \in \timeSet \label{eq:chapter3:trans3}\\
& \sum_{\unit \in \unitSet_n} \production_{\unit \timeUnit} + \sum_{r \in R_n} p_{rt} + \sum_{s \in S_n}  p_{st} + \nodeInjection_{\node \timeUnit} = \demand_{\node \timeUnit} , \timeUnit \in \timeSet, \node \in \nodeSet \label{eq:chapter3:nodalMarketClearing} \\
&\unitCommit_{gt} - \unitCommit_{gt-1} = \decisionStart_{gt} - \decisionStop_{gt} , \timeUnit \in \timeSet, \unit \in \unitSet \label{eq:chapter3:logic}\\
&  u_{gt},v_{gt},w_{gt} \in \{0,1\}, p_{gt},p_{rt},p_{st},pe_{st},inj_{nt},f_{lt} \in \mathbb{R} \label{eq:chapter3:domain}
\end{align}
\end{subequations}
(\ref{eq:chapter3:objective}) is the objective function of the UC consisting of the generation cost and start cost, $a_g$ is the constant cost, $b_g$ is the linear cost coefficient and $c_g$ is the quadratic cost coefficient. Constraint (\ref{eq:chapter3:pmin}) and (\ref{eq:chapter3:pmax}) ensure the minimum and maximum production of generators. Constraint (\ref{eq:chapter3:MDT1}) and (\ref{eq:chapter3:MDT2}) ensure the minimum up and downtime of generators. Constraint (\ref{eq:chapter3:tightRamp1})  and (\ref{eq:chapter3:tightRamp2}) ensure the ramping limits of generators between time steps. Constraint (\ref{eq:chapter3:RES}) ensures that the RES production is lower than the availability at that hour. (\ref{eq:chapter3:ChargeLimit}), (\ref{eq:chapter3:DischargeLimit}) and (\ref{eq:chapter3:EnergyLimit}) ensure the charge, discharge and energy storage limits for storage units. Equation (\ref{eq:chapter3:StorageProduction}) is the sum of charge and discharge i.e. the net storage production. Equation  (\ref{eq:chapter3:NetStorageProduction}) describes the relation between the charge, discharge and net power production of a storage unit.
 Equation (\ref{eq:chapter3:logic}) describes the logic between the binary commitment, start and stop variables of the generators.
 Equation (\ref{eq:chapter3:trans1}) describes the relation between the flow on transmission lines and the power injection at nodes. Constraint (\ref{eq:chapter3:trans3}) ensures flow limits on transmission lines.
 Equation (\ref{eq:chapter3:nodalMarketClearing}) ensures that the total generation meets the total demand at every node and time step.
At last, the commitment variables are binary while the generation are real numbers (\ref{eq:chapter3:domain})

\section{ADMM Example on simple UC problem}
To illustrate the application of ADMM to UC, we present here a simple UC example. Suppose we want to find the least cost power production schedule for a power system with a demand $D$ and $n$ generators with generation levels $p_i$  and a convex generation cost function $f_i$ for generator $i$.  This optimization problem can be formulated as:

\begin{align}
    \min \sum_{1 \leq i \leq n} f_i(p_i) \label{eq:chapter3:basicObjective}\\
    s.t. \nonumber\\
    D- \sum_{1 \leq i \leq n} p_i = 0 \label{eq:chapter3:basicConstraint}
\end{align}

We can solve this problem by first taking the augmented Lagrange relaxation of the demand coupling constraint (\ref{eq:chapter3:basicConstraint}) and for a giving $\rho$ turning it in the following function:

\begin{align}
\label{eq:chapter3:SmallLangragian}
    \mathcal{L}(p_1,\dots , p_n,\lambda) = \sum_{1 \leq i \leq n} f_i(p_i)+ \lambda ( D - \sum_{1 \leq i \leq n} p_i ) + \frac{\rho}{2}( D - \sum_{1 \leq i \leq n} p_i)^2
\end{align}

we can define the corresponding dual problem as:

\begin{align} \label{eq:chapter3:dualFunction}
&  \max_\lambda g(\lambda) \\
&  \text{where } g(\lambda) = \inf_{\mathbf{p}} \mathcal{L}(\mathbf{p},\lambda) \nonumber
\end{align}

If strong duality holds,  the optimal value of the dual problem (\ref{eq:chapter3:dualFunction}) is the same as the optimal value of the primal problem (\ref{eq:chapter3:basicObjective}). Then the optimal primal solution can be recovered from the optimal dual solution. We can find the optimal multiplier $\lambda$ by doing dual descent, i.e. iteratively solving $\inf_{\mathbf{p}} \mathcal{L}(\mathbf{p},\lambda)$ for a given $\lambda$ and update $\lambda$ in the direction of the gradient of $g$ (\ref{eq:chapter3:dualFunction}).

However, contrary to solving a standard Lagrangian iteration, we cannot determine the optimal generation $p_i$ independently of the other production variables  when solving an augmented Lagrangian iteration. Unfortunately,  this makes the augmented Lagrangian dual iteration as difficult as solving the original primal problem.

A solution to deal with the interdependence of the subproblems is ADMM.  The idea of ADMM is that you iteratively optimize one subproblem  while the other subproblems stay fixed. In other words, only a subset of the decision variables can change while the rest of the decision variables remain fixed  at the value of the previous iteration or, if they have already been updated, of the current iteration.  At the end of the ADMM iteration the dual multipliers are updated (see \cite{boyd2011distributed} for a detailed explanation of ADMM).

When the original problem is convex and the problem is split into two sets of variables then this method will converge to a primal feasible solution with an optimal value \cite{boyd2011distributed}. However, a lot of optimization problems have a decomposable structure that consist of multiple sets or blocks. To decompose a problem that has a multi-block structure we can either use multi-block Gauss Seidel ADMM or variable splitting ADMM which at its core has a two-block structure \cite{liu2015multi}. The former ADMM method has been applied to practical problems \cite{liu2015multi}  but lacks the converging guaranty of the latter \cite{chen2016direct}.

\subsection{Multi-block Gauss Seidel}\label{sect:GauseSeidelADMM}

We apply Gauss Seidel ADMM to determine the optimal Lagrange multipliers and generation for the UC example. At each iteration $k$ we sequentially minimize $p_i$ while keeping the other production values constant. The values $p_{j}^{k+1}$ are used for generators $j$ where $j<i$  and for generators $j>i$ we use the old value $p_j^{k}$ of the previous iteration, resulting in each iteration in $n$ statements of the form:
\begin{align}
      p_i^{k+1} \leftarrow \text{argmin}_{p_i}  \mathcal{L}(p_1^{k+1},\dots,p_{i-1}^{k+1},p_i,p_{i+1}^{k},\dots,p_n^k,\lambda^k) \label{eq:chapter3:productionOptimization}
\end{align}
The strength of ADMM is that we remove the interdependence created by the cross term in  (\ref{eq:chapter3:SmallLangragian}) by optimizing only for $p_i$ while keeping the other values constant.

At the end of each iteration, the multiplier is updated in the direction of the gradient of the relaxed constraint violation with step size $\rho$. For readability, we  define the relevant residual load $R_{i}$ for generator $i$:
\begin{align}
    R_i =  D  -  \sum_{1\leq j<i} p_{j}^{k+1}  -  \sum_{i<j\leq n} p_{j}^{k}
\end{align}
The full Gauss-Seidel multi-block ADMM that solves our small UC  example with $n$ generators is presented in Algorithm~\ref{alg:ADMMGaussSeidel}.

\begin{algorithm}[!htbp]
\begin{algorithmic}
\While{Stopping criteria have not been met}
  \State  $ p_1^{k+1}  \leftarrow  \text{argmin}_{p_1} \{f_1(p_1) - \lambda^k  p_1  + \frac{\rho}{2}( R_{1}  -  p_1 )^2 \} $ \Comment{$ \text{argmin}_{p_1}  \mathcal{L}(p_1,\dots,p_n^k,\lambda) $}
  \State $\vdots $
\State   $ p_i^{k+1} \leftarrow \text{argmin}_{p_i}  \{ f_i(p_i) - \lambda^k  p_i+ \frac{\rho}{2}(  R_i  -  p_i )^2 \}$ \Comment{$ \text{argmin}_{p_i}  \mathcal{L}(p_1^{k+1}\dots p_i\dots p_n^k,\lambda) $}
 \State  $ \vdots $
\State $   p_n^{k+1}  \leftarrow  \text{argmin}_{p_n}  \{ f_n(p_n) - \lambda^k  p_n + \frac{\rho}{2}(  R_n  -  p_n )^2 \} $ \Comment{$ \text{argmin}_{p_n}  \mathcal{L}(p_1^{k+1},\dots,p_n,\lambda) $}
\State $   \lambda^{k+1} \leftarrow   \lambda^{k} + \rho  ( D - \sum_{1\leq i \leq n} p_{i}^{k+1}) $
 \State $k \leftarrow k +1$
\EndWhile
\end{algorithmic}
\caption{Gauss-Seidel Multi-block ADMM}
\label{alg:ADMMGaussSeidel}
\end{algorithm}

An alternative method of applying ADMM on a problem with a multi-block structure is by variable splitting, but we found that this ADMM performed much worse than the Gauss Seidel ADMM. More information is presented in the appendix.

\section{New ADMM Algorithm for UC}
After the application of ADMM to a simple UC example, we now apply ADMM to the full UC problem described in Section 2. Although the procedure is also based on the relaxation of the coupling constraint that links generation and demand, the individual subproblems are more complicated.

In the UC description (\ref{eq:chapter3:objective}) - (\ref{eq:chapter3:domain}) the only constraint that links the generation and demand is (\ref{eq:chapter3:nodalMarketClearing}). We obtain the following augmented Lagrangian by relaxing this constraint:
\begin{subequations}
\begin{align}
    &     \mathcal{L}(\overrightarrow{p_{gt}},\overrightarrow{p_{rt}},\overrightarrow{p_{st}},\overrightarrow{inj_{nt}}, \overrightarrow{\lambda_{nt}}) =  \sum_{t \in T}\sum_{g \in G} \left[ a_g + b_g p_{gt} + c_g p_{gt}^2 + v_{gt} cost_{start}\right] + \nonumber\\ 
    & \sum_{t \in T} \sum_{n\in N} \left[ \lambda_{nt} RD_{nt}  +\frac{\rho}{2} RD_{nt}^2 \right]\\
    & \text{s.t. } (\ref{eq:chapter3:pmin}) - (\ref{eq:chapter3:trans3}) 
\end{align}
\end{subequations}
where:
\begin{align}
        RD_{nt} &=     \demand_{\node \timeUnit}  -  (    \sum_{\unit \in \unitSet_n} \production_{\unit \timeUnit} +  \sum_{r \in R_n} p_{rt} +   \sum_{s \in S_n} p_{st} + \nodeInjection_{\node \timeUnit}  \nonumber ) 
        \end{align}
$\lambda_{nt}$ is the Lagrange multiplier of coupling  constraint of node $n$ at time $t$ and $\rho$ is the quadratic penalty coefficient on the residual demand.

\begin{figure}[!htp]
    \centering
    \includegraphics[scale=0.8]{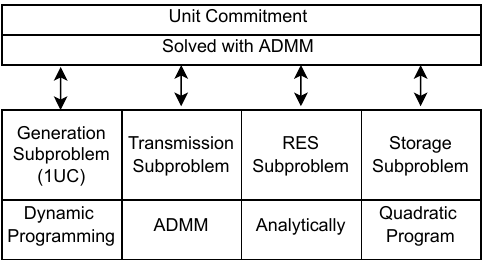}
    \caption{Flowchart of solving UC with ADMM by decomposing it into subproblems.}
    \label{fig:flowchart}
\end{figure}

At each iteration the optimal value of the decision variables will be determined while the rest of the decision variables remain fixed at the value of the previous iteration or, if they have already been updated, of the current iteration.
 Consequently, we iteratively solve the following subproblems: single unit commitment subproblem (1UC) for each generator $g \in G$, RES subproblem for each renewable resource $r \in R$, storage subproblem for each storage unit $s \in S$, and the transmission subproblem for each time step $t \in T$.
This is depicted in (Figure~\ref{fig:flowchart}).
After each iteration the Lagrange multipliers are updated in the direction of the sub-gradient which equals to the residual demand $RD_{nt}$. The ADMM Gauss-Seidel iteration proceeds according to Algorithm~\ref{alg:ADMMIter}.

\begin{algorithm}[!htp]
\centering
\caption{Gauss-Seidel ADMM iteration}
\label{alg:ADMMIter}
\begin{algorithmic}
\ForAll{$g \in G$}
\State $\mathbf{p}_g \leftarrow \text{argmin}_{p_g} \mathcal{L}(\dots,p_{g-1}^{k+1},p_{g},p_{g+1}^{k}\dots, \lambda^{k})$  \Comment{Solve 1UC subproblem}
\EndFor
\ForAll{$r \in R$}
\State $\mathbf{p}_r \leftarrow \text{argmin}_{p_r} \mathcal{L}(\dots,p_{r-1}^{k+1},p_{r},p_{r+1}^{k}\dots, \lambda^{k}) $  \Comment{Solve 1RES subproblem}
\EndFor
\ForAll{$s \in S$}
\State $p_s \leftarrow \text{argmin}_{p_s} \mathcal{L}(\dots,p_{s-1}^{k+1}p_s,p_{s+1}^{k}\dots, \lambda^{k}) $  \Comment{Solve 1Storage subproblem}
\EndFor
\ForAll{$t \in T$}
\State $\mathbf{inj}_t \leftarrow \text{argmin}_{inj_t }\mathcal{L}(\dots,inj_t,\dots, \lambda^{k})$   \Comment{Solve Transmission subproblem}
\EndFor
\ForAll{$t \in T$}
\ForAll{$n \in N$}
 \State $ \lambda_{nt}^{k+1} \leftarrow \lambda_{nt}^{k} + \rho RD_{nt}^{k}$ \Comment{Update Lagrangian  multipliers}
\EndFor
\EndFor 
\end{algorithmic}
\end{algorithm}

\newcommand{\oneUC}[0]{\text{1UC}_g(\boldsymbol{\lambda},R_{p_{gt}})}
\newcommand{\oneRES}[0]{\text{1RES}_r(\boldsymbol{\lambda},R_{rnt})}
\newcommand{\oneStorage}[0]{\text{1Storage}_s(\boldsymbol{\lambda},R_{p_{st}})}
\newcommand{\oneTrans}[0]{\text{Trans}_s(\boldsymbol{\lambda},R_{inj_{nt}})}

\subsection{Solving the subproblems}

\subsubsection{Generator subproblem 1UC}

The 1UC problem is the problem of determining the optimal power schedule for a single  generator.
In each iteration of our Gauss-Seidel algorithm, we have to solve a number of 1UC problems (see Line 2 in Algorithm~\ref{alg:ADMMIter}). For a given generator $g$ at node $n$, this problem can be defined as:
\begin{align}\label{eq:chapter3:1UCLagrange}
    \text{argmin}_{p_{g}} \mathcal{L}(\dots,p_{g-1}^{k+1},p_{g},p_{g+1}^{k}\dots, \lambda^{k}) \text{ s.t. } (\ref{eq:chapter3:pmin}) - (\ref{eq:chapter3:domain})
    \end{align}
Only $p_g$ is variable here, the other production variables $p_g$ $(g \in G_n)$ are constant.
 
Similar to the the small UC example in Section~\ref{sect:GauseSeidelADMM}, we define the relevant residual load $R_{gnt}$  for determining the optimal power output of generator $g$ at node $n$ and time $t$ by
\begin{align}
    R_{gnt} = \demand_{\node \timeUnit}  -  \sum_{1\leq i<g} \production_{i t}^{k+1}  -  \sum_{g<i\leq  \mid G_n \mid } \production_{i t}^{k} -  \sum_{r \in R_n} p_{rt}^k    - \sum_{\sUnit \in \sUnitSet_\node} p_{\sUnit \timeUnit}^k  - \nodeInjection_{\node \timeUnit}^k
\end{align}
We can rewrite (\ref{eq:chapter3:1UCLagrange}) to an equivalent problem formulation:
\begin{align}
& \label{eq:chapter3:1UCsub}\min \sum_{t \in T} a_g + b_g  p_{gt} + c_g p_{gt}^2 + v_{gt} cost_{start} - \lambda_{t} p_{gt}   + \frac{\rho}{2} (R_{gnt} -  p_{gt})^ 2\\
& \text{subject to } (\ref{eq:chapter3:pmin}),(\ref{eq:chapter3:pmax}) ,(\ref{eq:chapter3:MDT1}),(\ref{eq:chapter3:MDT2}) ,(\ref{eq:chapter3:logic}) ,(\ref{eq:chapter3:tightRamp1}) ,(\ref{eq:chapter3:tightRamp2})
\end{align}

\paragraph{Solving 1UC.}

Since we have to solve a series of 1UC problems in each iteration, solving these problems very efficiently is crucial. Here we use our dynamic programming algorithm $RRF+$ \cite{wuijts2021improved}, which significantly outperforms earlier 1UC algorithms in terms of computation time. This advanced dynamic programming algorithm models part of the states by convex functions instead of single values. This algorithm is based on a recurrence relation on  functions that represents for each generator state at each time step the value of the optimal generator schedule that ends in that state. The RFF+ algorithm is very efficient because it does not need to compute the optimal economic dispatch for each possible ‘on-interval’ of a generator and it can identify superfluous functions which reduces the computation time significantly (see \cite{wuijts2021improved} for more details).

\subsubsection{Renewable generation subproblem 1RES}
 The 1RES subproblem is the problem of finding the optimal renewable energy generation given its availability and is given by:
\begin{align}
     &   \text{argmin}_{p_r} \mathcal{L}(\dots,p_{r-1}^{k+1},p_{r},p_{r+1}^{k}\dots, \lambda^{k})  &\text{ s.t. } (\ref{eq:chapter3:pmin}) - (\ref{eq:chapter3:domain})\label{eq:chapter3:1RESLagrange}
\end{align}
We define the relevant residual load $R_{rnt}$ for RES $r$
at node $n$ and time $t$.
\begin{align}
    R_{rnt} = \demand_{\node \timeUnit}  -  \sum_{g \in G_n}  \production_{g t}^{k+1}   -  \sum_{1\leq i<r} \production_{r t}^{k+1}  -  \sum_{r<i\leq  \mid R_n \mid } \production_{r t}^{k}  - \sum_{s \in S_n} p_{st}^k - \nodeInjection_{\node \timeUnit}^k
\end{align}
Now (\ref{eq:chapter3:1RESLagrange}) is equivalent to the following program:
\begin{align}
&\min \sum_{t \in T} -\lambda_{t}^k p_{rt}  + \frac{\rho}{2}  (R_{rnt} -  p_{rt})^ 2 \label{eq:chapter3:1RESsub}\\
&  \text{subject to }    p_{rt} \leq AF_{rt} \productionMax_{rt} &  t \in T 
\end{align}
We can easily solve this problem analytically since it is a simple minimization of a quadratic function with bounds, which results in:

\begin{equation}
 p^{k+1}_{rt} = 
    \begin{cases} 
      0 & \frac{\lambda^k_{t}}{R_{rnt} }< 0\\
   \frac{\lambda^k_{t}}{R_{rnt} }  &  0 \leq \frac{\lambda^k_{t}}{R_{rnt} } \leq AF_{rt} \productionMax_{rt}   \\
      AF_{rt} \productionMax_{rt}  & \frac{\lambda^k_{t}}{R_{rnt}} > AF_{rt} \productionMax_{rt} \\
   \end{cases}
\end{equation}

\subsubsection{Storage subproblem 1Storage}
The storage subproblem to find the least cost charging and discharging strategy is given by:
\begin{align}
     &       \text{argmin}_{p_s} \mathcal{L}(\dots,p_{s-1}^{k+1},p_s,p_{s+1}^{k}\dots, \lambda^{k})  &\text{ s.t. } (\ref{eq:chapter3:pmin}) - (\ref{eq:chapter3:domain}) \label{eq:chapter3:1StorageLagrange}
\end{align}
We define the relevant residual load $R_{snt}$ for storage unit $s$ at node $n$ and time $t$.
\begin{align}
    R_{snt} = \demand_{\node \timeUnit}  -  \sum_{g \in G_n}  \production_{g t}^{k+1}  -  \sum_{r \in R_n} p_{rt}^{k+1}  -  \sum_{1\leq i<s}  p_{st}^{k+1} - \sum_{s<i\leq  \mid S_n  \mid } p_{st}^{k}- \nodeInjection_{\node \timeUnit}^k  
\end{align}
To find the best storage schedule (\ref{eq:chapter3:1StorageLagrange}), we can simply solve the following quadratic program:
\begin{align}
&\min \sum_{t \in T}- \lambda^k_{t} p_{st}   + \frac{\rho}{2} (R_{snt} - p_{st} )^2 \label{eq:chapter3:1STORsub}\\
& \text{subject to } (\ref{eq:chapter3:EnergyLimit}),(\ref{eq:chapter3:ChargeLimit}) ,(\ref{eq:chapter3:DischargeLimit}),(\ref{eq:chapter3:NetStorageProduction}),(\ref{eq:chapter3:StorageProduction})
\end{align}

Although other algorithms exist that solve this problem efficiently with dynamic programming \cite{bannister1991rapid}, we solve this storage subproblem with the QP-solver Gurobi.

\subsubsection{Transmission subproblem}

The transmission subproblem is given by
\begin{align}
     &   \text{argmin}_{inj_{nt} }\mathcal{L}(\dots,inj_{nt},\dots, \lambda^k) &\text{ s.t. } (\ref{eq:chapter3:pmin}) - (\ref{eq:chapter3:domain})\label{eq:chapter3:1TransLagrange}
\end{align}
The transmission subproblem (\ref{eq:chapter3:1TransLagrange}) is harder than the RES or storage subproblems but fortunately it consists of a set of the subproblems which are independent in time. We define the relevant residual load $R_{nt}$ for the optimal nodal injection at node $n$ and time $t$.
\begin{align}
    R_{nt} = \demand_{\node \timeUnit}  -  \sum_{g \in G_n}  \production_{g t}^{k+1}   - \sum_{r \in R_n} \production_{r t}^{k+1}  - \sum_{s \in S_n} p_{st}^{k+1} 
\end{align}
Moreover, we define a set $L_n$ as $\{l\mid l \in L, n' \in N, l = (n' \rightarrow n) \}$ i.e. the set of transmission line going into $n$. For each time step $t \in T$ we have to solve the transmission subproblem that is equivalent to  (\ref{eq:chapter3:1TransLagrange}):
\begin{subequations}
\label{Tsub}
\begin{align}
&\min \sum_{n \in N}  - \lambda^k_{n} inj_{n}   + \frac{\rho}{2} (R_{nt} - inj_{n}) ^ 2\\
& \text{subject to} \nonumber \\
& \nodeInjection_{\node  } = \sum_{l \in  L_{n}} \flow_{l} & \forall \node \in \nodeSet \label{trans1ForT} \\
&\underline{\flow_{l}} \leq \flow_{l} \leq \overline{\flow_{l}}  &  \forall l \in L \label{trans3ForT}
\end{align}
\end{subequations}
This subproblem is a quadratic program  that can be solved by a QP-Solver. However, computation time is long, since it needs to be solved for every time step at every ADMM iteration. We found that solving this subproblem with our own ADMM algorithm is much faster. Since (\ref{Tsub}) is a convex quadratic program we can solve it to optimality with a standard ADMM procedure by relaxing the injection and flow coupling constraint (\ref{trans1ForT}). This results in the following augmented Lagrangian:
\begin{align}
\label{TsubAugDual}
& \mathcal{L}^{trans}(\overrightarrow{\pi_n},\overrightarrow{inj_n},\overrightarrow{f_l})= \\
&\sum_{n \in N} \left(- \lambda^k_{n} inj_{n}  + \frac{\rho}{2} (R_{nt} - inj_{n}) ^ 2 + \pi_{n} (\nodeInjection_{\node  } - \sum_{l \in  L_{n}} \flow_l) + \frac{\rho_{trans}}{2} (\nodeInjection_{\node  } - \sum_{l \in L_n} \flow_l)^2 \right) \nonumber
\end{align}

When we solve (\ref{TsubAugDual}) with ADMM we iteratively need to solve two kinds of subproblems: determining the optimal injections and determining the optimal flow on transmission lines. The former is solved by finding the minimum of a quadratic function and the latter by finding the minimum of a quadratic function with bounds. 

\begin{algorithm}[!htbp]
\label{alg:ADMMTrans}
\caption{Transmission Subproblem  Gauss-Seidel ADMM}
\begin{algorithmic}[1]
\ForAll{$n \in N$}
\State $inj_{n} \leftarrow \text{argmin}_{\nodeInjection_{n}}  \mathcal{L}^{trans}(\dots,inj_{n-1}^{k+1},inj_n,inj_{n+1}^{k},\dots)$ \Comment{see (\ref{eq:chapter3:injectionSubproblem})}
\EndFor
\ForAll{$l \in L$}
\State $f_l \leftarrow \text{argmin}_{f_l}\mathcal{L}^{trans}(\dots,f_{l-1}^{k+1},f_l,f_{l+1}^{k},\dots)$\Comment{see (\ref{eq:chapter3:lineSubproblem})}
\EndFor
\ForAll{$n \in N$}
 \State $ \pi_{n} \leftarrow \pi_{n} + \rho_{trans} (\nodeInjection_{\node  } - \sum_{l \in  L_{n}} \flow_{l} )$ \Comment{Update Lagrangian  multipliers}
\EndFor
\end{algorithmic}
\end{algorithm}

When determining the optimal nodal injection of a single node $n$ at ADMM iteration $k$, the other nodal injection and flows are constant. Determining the optimal injection is done by finding the minimal point of the following quadratic function:
\begin{subequations}
\begin{align}\label{eq:chapter3:injectionSubproblem}
&\min \underbrace{( - \lambda^k_{n} - \rho R_{nt}+\pi_{n} -   \rho_{trans} \sum_{l \in L_n} \flow_{l}^k  )}_{\text{linear coefficient}} \textcolor{green}{inj_{n}} + \underbrace{( \frac{\rho}{2} + \frac{\rho_{trans}}{2})}_{\text{quadratic coefficient}} \textcolor{green}{inj_{n}} ^2
\end{align}
\end{subequations}

Similarly, determining the optimal flow for a single transmission line $l = (n',n)$ at the ADMM iteration, the other flows and nodal injections are constant. Determining the optimal flow is done by finding the minimal point of the following quadratic function within certain flow limits:

\begin{subequations}
\begin{align}\label{eq:chapter3:lineSubproblem}
&\underbrace{(-\pi_{n'} + \pi_{n} - \rho_{trans} (\nodeInjection_{n'}^{k+1} - \sum_{l' \in L_{n'} \setminus \{ l\}} f_{l'}^k ) + \rho_{trans} (\nodeInjection_{n}^{k+1} - \sum_{l' \in L_{n} \setminus \{ l\}} f_{l'}^k )}_{\text{linear coefficient}} \textcolor{green}{f_{l}} + \nonumber   \\
& \underbrace{\rho_{trans}}_{\text{quadratic coefficient}} \textcolor{green}{f_l}^2 \\
& \text{subject to} \nonumber \\
&\underline{f_l} \leq f_l \leq \overline{f_l}
\end{align}
\end{subequations}

\subsection{Our ADMM algorithm}
When all the subproblems are convex,  ADMM will converge to the optimal solution after some  iterations\footnote{For the optimal Exchange this is guaranteed and for the Gauss Seidel variant this is probable \cite{liu2015multi}, i.e. there are theoretical cases where it does not converge but experimentally there have been success. In our experiments with multiple convex UC instances the Gauss Seidel variant always converged to the optimal solution.}. But since 1UC is non-convex, optimality or convergence are not guaranteed.  This means that after a certain number of iterations the residual load is still significant indicating that the current solution is infeasible due to the violation of the demand coupling constraint.

\begin{figure}[!htbp]
    \centering
        \includegraphics[scale=0.8]{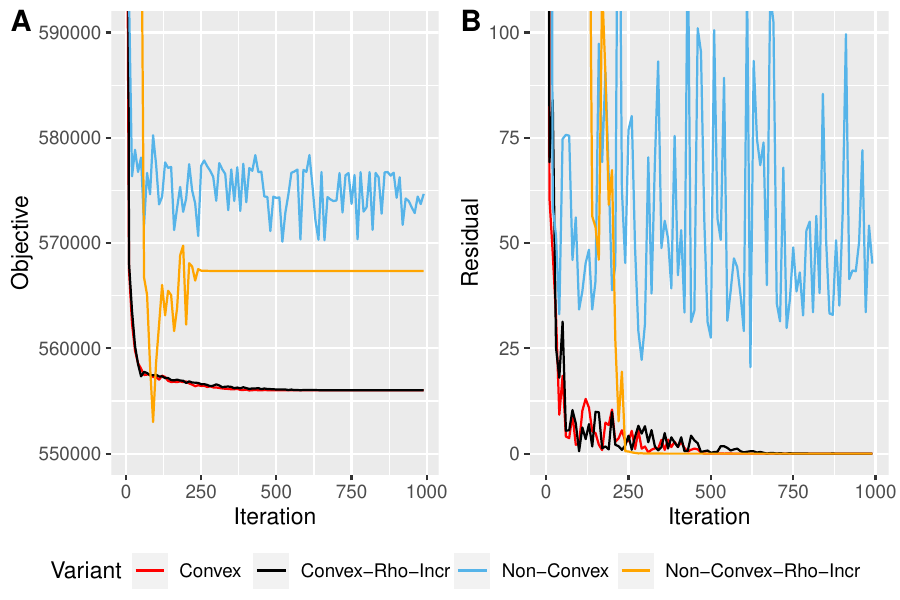}
    \caption{1000 ADMM iterations on a convex and a non-convex UC instance with and without increasing $\rho$.}
    \label{fig:converge}
\end{figure}

To solve the nonconvergence,  we iteratively increase the penalty coefficient $\rho$ to force convergence and feasibility. In \cite{xu2017no}, it was shown that this helps to achieve faster convergence of convex problems.  Figure~\ref{fig:converge} shows in red the convergence of the solution when 1UC is convex\footnote{We convexified the 1UC problem by removing the binary requirement on the commitment variables.} and the nonconvergence (at least not after 1000 iterations) in blue of the non-convex problem. It also shows the convergence of the non-convex problem in orange when $\rho$ is increased.
As $\rho$ grows to $\infty$ the constraint violation is encouraged to go to $0$. We start with a small $\rho$ and then iteratively increase it after $m$ iterations. We found that the slower this process goes the better the heuristic solution becomes. The whole procedure is presented in Algorithm~\ref{alg:ADDM}.

\begin{algorithm}[!htbp]
\caption{ADMM Increasing $\rho$}

\label{alg:ADDM}
\begin{algorithmic}[1]
\State $\rho \leftarrow$ initial $\rho$ 
\While{  $\sum_{n\in N, t\in T} RD_{nt}> \epsilon$ $\land$  $k < $Max iterations}
\For{$i \in \{1,\dots,m\}$}
  \State Gauss Seidel ADMM Iteration
\State k $\leftarrow$ k + 1
\EndFor
\State $\rho \leftarrow \alpha \rho $
\EndWhile
\end{algorithmic}
\end{algorithm}

Here $m$ is the number of iterations after we increase $\rho$, by multiplying $\rho$  times $\alpha$. Varying $\alpha$ and $m$ gives us a way to increase and decrease the rate of convergence and therefore the total computation time but at the expense of solution quality.

 
The intuition is as follows, when solving the subproblems the individual generators, RES and storage units are encouraged to produce power at every time step by the Lagrangian multipliers
(when positive) and encouraged to minimize the total residual by the quadratic penalty. Our ADMM algorithm could be interpreted as a greedy algorithm that first steers to the (sub)optimal solution by the Lagrangian multipliers and then steers towards feasibility by increasing   $\rho$. 
Theoretically not much is known about converges or the quality of the solution when applying ADMM to non-convex problems \cite{boyd2011distributed}.  In the next section we will show experimentally that our proposed algorithm gives a fast  feasible solution of good quality and scales well with the number of times steps.


\section{Computational Experiments}

\begin{table}[!htbp]
\small\centering
\begin{tabular}{l|l|l|l|l|l}
       Name   & \rot{Qudractic Cost} & \rot{Transmission} & \rot{Storage} & \rot{RES} & \rot{Source} \\\hline
GA10      & \Yes       & \No  & \No & \No & \cite{kazarlis1996genetic}\cite{rahman2014metaheuristic}      \\
TAI38     & \Yes       & \No  & \No & \No & \cite{huang1997new}                                           \\
RCUC50    & \Yes       & \No  & \No & \No & \cite{frangioni2009tighter}\\
GMLC73    & \No        & \No  & \No & \Yes& \cite{barrows2019ieee}\\
A110      & \Yes       & \No  & \No & \No & \cite{orero1997large}                                         \\
KOR140    & \Yes       & \No  & \No & \No & \cite{park2010improved}\cite{moradi2015semi}                  \\
OSTRO187  & \Yes       & \No  & \No & \No & \cite{ostrowski2012tight}                                  \\
RCUC200   & \Yes       & \No  & \No & \No & \cite{frangioni2009tighter} \\
HUB223    & \No        & \No  & \No & \Yes& \cite{huber2015modeling}\\
CA426     & \No        & \No  & \No & \Yes& \cite{knueven2018mixed}\\
FERC923   & \No        & \No  & \No & \Yes& \cite{knueven2018mixed}\cite{krall2012rto}\\
RTS26     & \Yes       & \Yes & \No & \No & \cite{wang1995short}\cite{wang1993effects}\\
RTS54     & \No        & \Yes & \No & \Yes& \cite{huber2015modeling}\\
RTS96     & \No        & \Yes & \No & \No & \cite{pandzicunit}\\
DSET304   & \No        & \Yes & \Yes& \Yes& \cite{kavvadias2018integrated}                                       
\end{tabular}
\caption{Overview of the instances in the benchmark
set. All instances have thermal generators with
techno-economic characteristics but only some have
a quadratic generation cost,  a transmission system, storage or a renewable
energy with availability factors. See \cite{wuijts2022modelchar} for more detail.}
\label{table:instances}
\end{table}
To test the efficiency of our proposed algorithms we have performed multiple experiments on 15 well-known benchmark instances\footnote{\url{https://github.com/rogierhans/UCBenchmark}} found in the literature (Table~\ref{table:instances}) created by \cite{wuijts2022modelchar}. We solved the UC instances with our Gauss Seidel ADMM algorithm and with Gurobi with the presented MIL(Q)P formulation (\ref{eq:chapter3:objective}) - (\ref{eq:chapter3:domain}).  Next, we compared the quality of the solutions and the computation times of these algorithms. 

 First, we experimented with different values of alpha, namely  $1.01$, $1.05$, $1.1$ or, $1.2$. Due to the large number of experiments we needed to perform we chose to omit the instances with transmission. Note that as $\alpha$ increases that the computation time and number of iterations decreases. Finally, we performed more detailed experiments for the average value $\alpha = 1.1$. For the last  experiments we also included the instances with transmission. 

All experiments were repeated 10 times to account for the random order of updating the subproblems and the random initial Lagrange multipliers. We set the time horizon to $24,48,\dots,168$ time steps. In all runs, $\rho$ was initially set at $0.0001$ and $m$ at $1$. If $\rho$ was set lower or $m$ higher then the solution quality is better but the algorithm takes more iterations to converge. In preliminary  experiments we found that these parameters give a good trade-off between computation time and solution quality except for 1 instance which will be addressed in the next section.  Both ADMM algorithms and Gurobi had a time limit of 1 hour.

\begin{figure}[!htbp]
    \centering
    \includegraphics[scale=0.35]{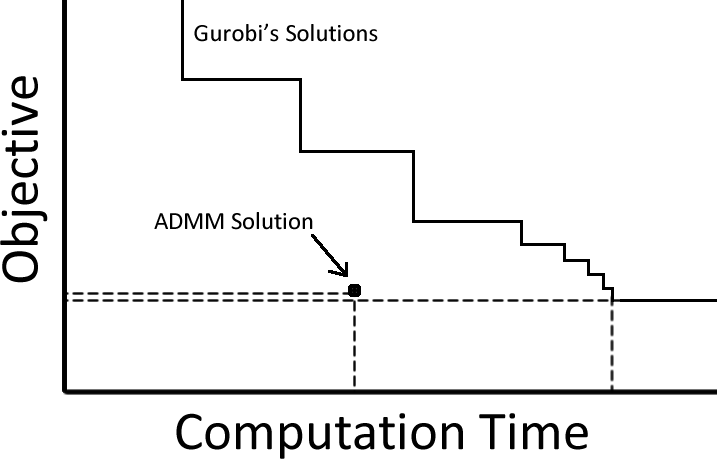}
    \caption{Comparing the ADMM solution to Gurobi.}
    \label{fig:GurobiCompare}
\end{figure}

Gurobi’s MIL(Q)P solver continuously finds heuristic upper bounds and lower bounds until the specified MIP-gap is reached. Our ADMM algorithm produces a single solution, which turned out to be feasible in all our experiments. To compare our ADMM algorithms to Gurobi, we ran the MIL(Q)P model with a MIP-gap of 0 and recorded the intermediate heuristic solutions along with a timestamp at which time they were found. We can then compare the computation times of our ADMM algorithm to the time that Gurobi needs to find a heuristic solution with at most the same value (see Figure~\ref{fig:GurobiCompare}).

\subsection{Results}

\begin{figure}[!htbp]
    \centering
    \includegraphics[scale=0.8]{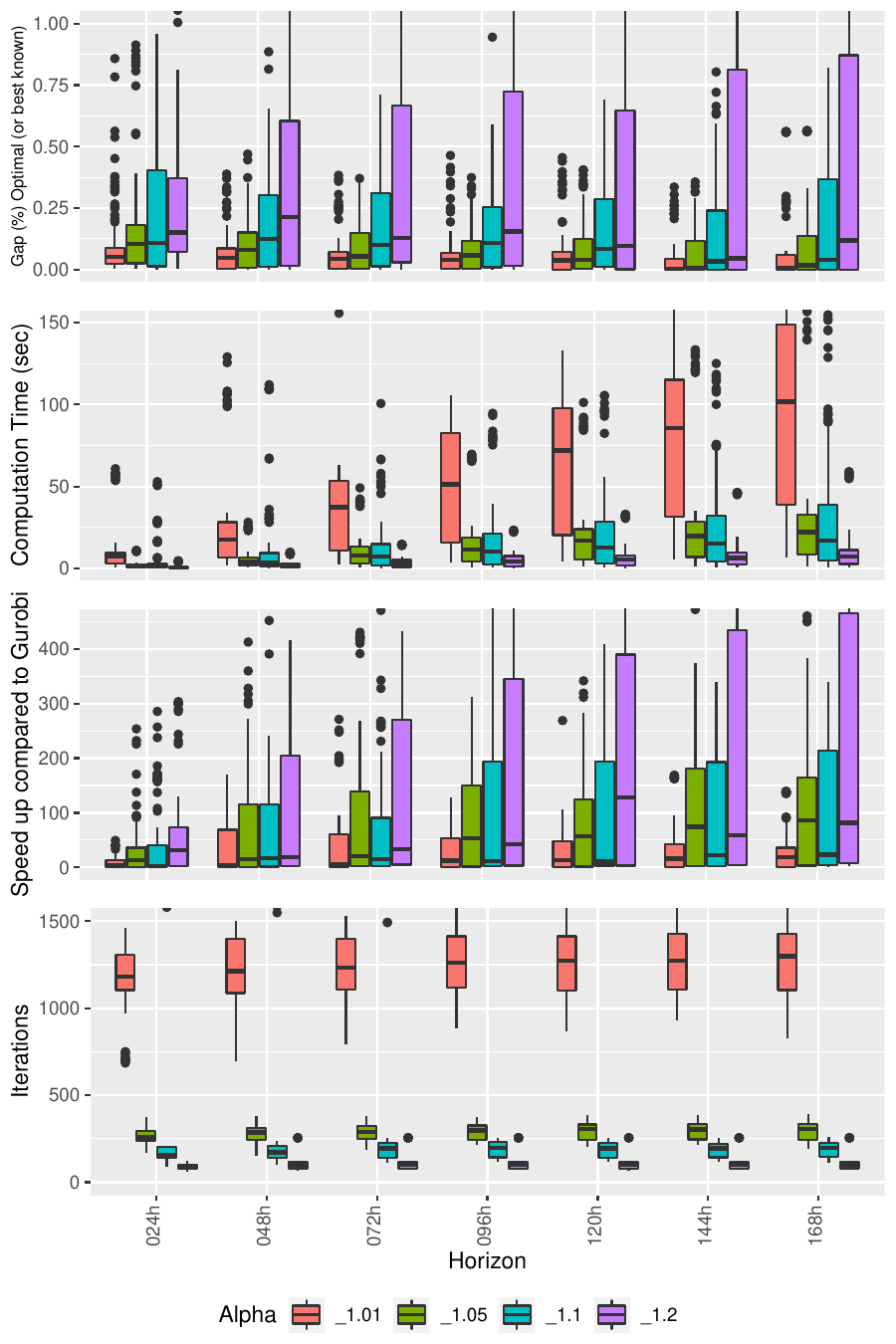}
    \caption{The gap, computation time, speed-up and the number of iterations for varying time horizons and values of $\alpha$ for the 11 instances without transmission.}
    \label{fig:summary}
\end{figure}

Figure~\ref{fig:summary} and Table~\ref{tab:Summary} show the aggregated results for the instance without transmission and for varying values of $\alpha$. We present the number of iterations, the gap with the optimal (or best known) solution, and the speed-up factor compared to solving the problem by Gurobi.  

The figure shows that with increasing $\alpha$ the computation time and number of iterations decreases but the solution quality also gets worse. An $\alpha$ of $1.01$ for example results in small gap indicating high quality solutions but it takes some time for $\rho$ to be of sufficient size that feasibility and convergence is forced. An  $\alpha$ of $1.2$ has the opposite effect, the  algorithm converges fast, for one instance it was 5801.3 times faster for Gurobi to find a solution of similar quality. However, the average gap is rather large and it has outliers, for example, one solution for UC instance \textit{KOR140} was 3.7\% removed from the optimal value. 

\begin{table}[!htbp]
\footnotesize
\begin{tabular}{c|cccc|cccc|cccc}
$\alpha$ & \multicolumn{4}{c|}{Iterations} & \multicolumn{4}{c|}{Gap \%} & \multicolumn{4}{c}{Speedup}    \\
         & avg    & median  & min  & max   & avg   & median & min & max  & avg   & median & min  & max    \\ \hline
1.01     & 1229   & 1239    & 685  & 1696  & 0.07  & 0.03   & 0   & 0.86 & 29.6  & 6.0    & 0.09 & 271.9  \\
1.05     & 287    & 292     & 154  & 393   & 0.09  & 0.04   & 0   & 0.91 & 108.4 & 19.4   & 0.5  & 1813.1 \\
1.1      & 163    & 168     & 90   & 216   & 0.15  & 0.08   & 0   & 1.15 & 195.3 & 35.6   & 0.8  & 3109.4 \\
1.2      & 104    & 100     & 60   & 255   & 0.54  & 0.13   & 0   & 3.7  & 310.4 & 43.3   & 1.0  & 5801.3
\end{tabular}
\caption{Summary of the experiments with varying values of $\alpha$ for the 11 instances without transmission.}
\label{tab:Summary}
\end{table}

\begin{figure}[!htbp]
    \centering
    \includegraphics[scale=0.9]{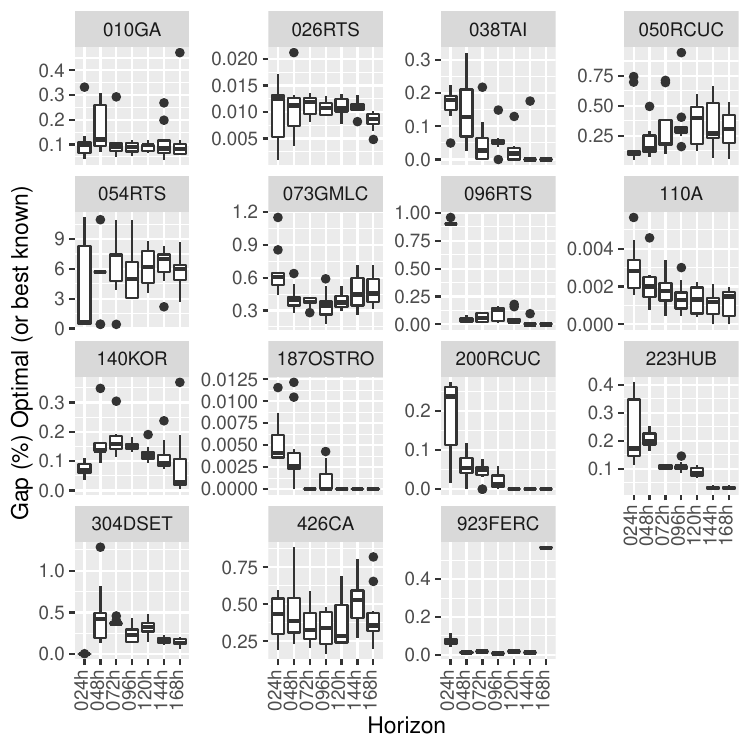}
    \caption{The gap of the objective value of the solution produced by our ADMM algorithm relative to the optimal or best known solution. For these experiments $\alpha$ was set at 1.1.}
    \label{fig:Ratio}
\end{figure}

\begin{figure}[!htbp]
    \centering
    \includegraphics[scale=0.9]{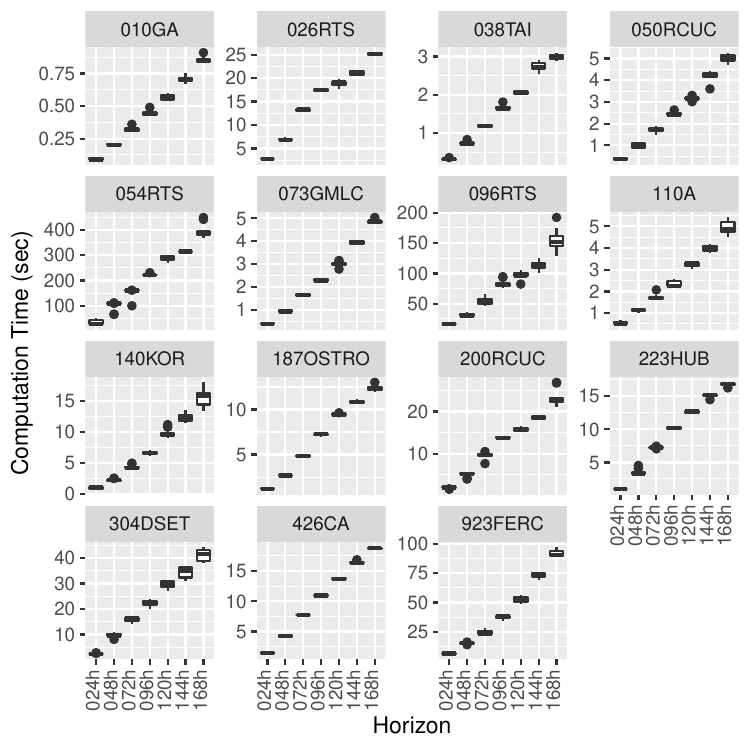}
    \caption{The computation time of all instances with a varying length of the time horizon. For these experiments $\alpha$ was set at 1.1 and for every instance and was repeated 10 times.} 
    \label{fig:Time}
\end{figure}

\begin{figure}[!htbp]
    \centering
    \includegraphics[scale=0.9]{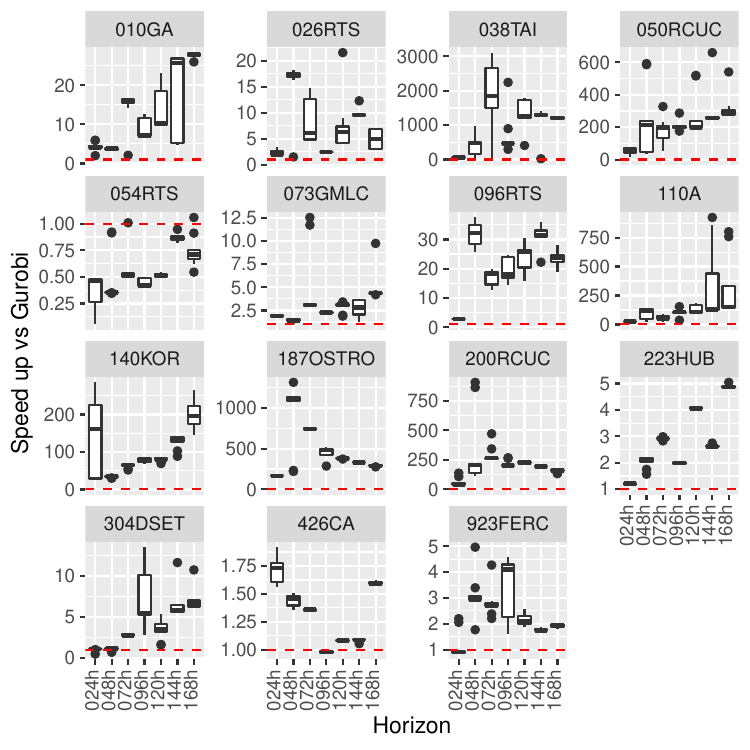}
    \caption{The speedup factor compared to the production of an equal solution by Gurobi.  For these experiments $\alpha$ was set at 1.1 and for every instance in was repeated 10 times.}
    \label{fig:SpeedUp}
\end{figure}

Figure~\ref{fig:Ratio}, Figure~\ref{fig:Time} and Figure~\ref{fig:SpeedUp} show the results of the experiments where $\alpha$ was 1.1 in more detail. These results also include the four instances with transmission.  Detailed results for $\alpha$ set at $1.01$, $1.05$ or $1.2$ would show similar results but with a faster or slower running time and better or worse solution quality according to the trend presented in Table~\ref{tab:Summary}. 

Figure~\ref{fig:Ratio} shows that for most instances the quality is good, for others it is a little worse but still close to optimal while for the instance \textit{RTS54} the solutions are far away from the optimum. Even with those outliers on average the algorithm performs well, on average 0.08\% away from the optimal value or the best known (Table~\ref{tab:Summary}). Moreover, if we choose a different value for $\alpha$  then these gaps would be a lot smaller.  In some instances,  the solutions produced by our algorithm was the best known. For the instances  \textit{TAI138}, \textit{RCUC200},  \textit{OSTRO178} and \textit{KOR140} Gurobi could not find a better solution within the 1 hour time limit for time horizons $>48$.

Figure~\ref{fig:Time} presents the computation time increase of our ADMM Gauss Seidel Algorithm when the time horizon increases. For most instances it seems to grow linearly with the number of time steps while for others the growth seems slightly more than linear. Recall from Figure~\ref{fig:summary} that the number of iterations needed for convergence stays constant when the time horizons increases. However, the computation time of the underlying 1UC subproblem grows when the time horizon increases\footnote{The time complexity of our 1UC algorithm was experimental established as linear \cite{wuijts2021improved} in the length of the time steps but also linear in the length of the minimum up and downtime. This is speculative but the superlinearity of \textit{FERC923} could be explained by the fact that this instance has multiple units with minimum up and downtime of 168 hours. As the time horizon increases the computation time of the subproblem increases quadratically up until the point that the minimum up- and downtime is achieved but this needs more investigation in order to be confirmed. }.

Figure~\ref{fig:SpeedUp} shows how fast our algorithm can find a high quality solution compared to Gurobi. Speedup is defined as the ratio of the computation time of our algorithm and of Gurobi to find a similar solution. The computation time of Gurobi is set at the first time it finds a solution with an objective equal or better than our final objective. For most instances this ratio is large and indicates a large speedup. It seems that our algorithm performs much better than Gurobi on instances with quadratic generation cost. For some instances, especially the ones with linear generation cost, the ratio is smaller and even sometimes below 1. These are for example, \textit{RTS54}, \textit{GMLC73},  \textit{HUB223}, \textit{DSET304} \textit{CA426},  and \textit{FERC923}. This highlights that for Gurobi solving a MILP is much faster that solving a MIQP. Still even for the instances with linear generation cost the speedup obtained by our algorithm is  significant. Moreover, for some cases the speedup increases when the time horizon increases.  

Because in the cases Gurobi could not find an equal or better solution within the time limit of 1 hour, we set the computation time to one hour, the real ratios for those instances are larger. Details are presented in Table \ref{tab:Gurobi}. This  also explains the seemingly downward trend of the instances \textit{RCUC50}, \textit{KOR140}, \textit{OSTRO187} and \textit{RCUC200}.

\begin{table}[!htbp]
\centering
\begin{tabular}{c|ccccccc}
       Time horizon     & 24h      & 48h      & 72h      & 96h      & 120h      & 148h      & 168h      \\ \hline
GA10        & \it{1s}      & 6s      & 35s     & 76s     & 171s     & 384s     & 551s     \\
RTS26       & 59s     & 0.006\% & 0.029\% & 0.069\% & 0.1\%    & 0.107\%  & 0.111\%  \\
TAI38       & 0.516\% & 1.356\% & 1.746\% & 1.822\% & 1.922\%  & 2.078\%  & 2.166\% \\
RCUC50      & 0.031\% & 0.106\% & 0.142\% & 0.177\% & 0.184\%  & 0.202\%  & 0.25\%   \\
RTS54       & 0.031\% & 0.084\% & 0.147\% & 0.21\%  & 0.236\%  & 0.296\%  & 0.328\%  \\
GMLC73    & 3s      & 12s     & 188s    & 0.044\% & 0.08\%   & 0.027\%  & 0.078\%  \\
RTS96       & 1.601\% & 1.766\% & 1.835\% & 1.89\%  & 1.921\%  & 2.062\%  & 2.536\%  \\
A110        & 0.001\% & 0.009\% & 0.012\% & 0.014\% & 0.015\%  & 0.017\%  & 0.02\%   \\
KOR140      & 0.64\%  & 0.929\% & 1.026\% & 1.06\%  & 1.109\%  & 1.152\%  & 1.206\%  \\
OSTRO187    & 0.029\% & 0.05\%  & 0.063\% & 0.07\%  & 0.115\%  & 0.114\%  & 0.337\%  \\
RCUC200     & 0.091\% & 0.163\% & 0.283\% & 0.308\% & 11.582\% & 20.613\% & 31.123\% \\
HUB223      & 3s      & 0.001\% & 0.001\% & 0.002\% & 0.002\%  & 0.001\%  & 0.002\%  \\
DSET304     & 3s      & 1664s   & 0\%     & 0\%     & 0.031\%  & 0.03\%   & 0.063\%  \\
CA426 & 6s      & 62s     & 118s    & 197s    & 406s     & 615s     & 971s     \\
FERC923  & 94s     & 96s     & 389s    & 378s    & 3601s    & 1043s    & 0.001\%  \\
\end{tabular}
\caption{ Computation time for Gurobi to find the optimal solution in seconds.  If Gurobi failed to find the optimal solution within 1 hour computation time limit then the MIP-gap is shown.}
\label{tab:Gurobi}
\end{table}

\begin{figure}[!htbp]
    \centering
    \includegraphics[scale=0.7]{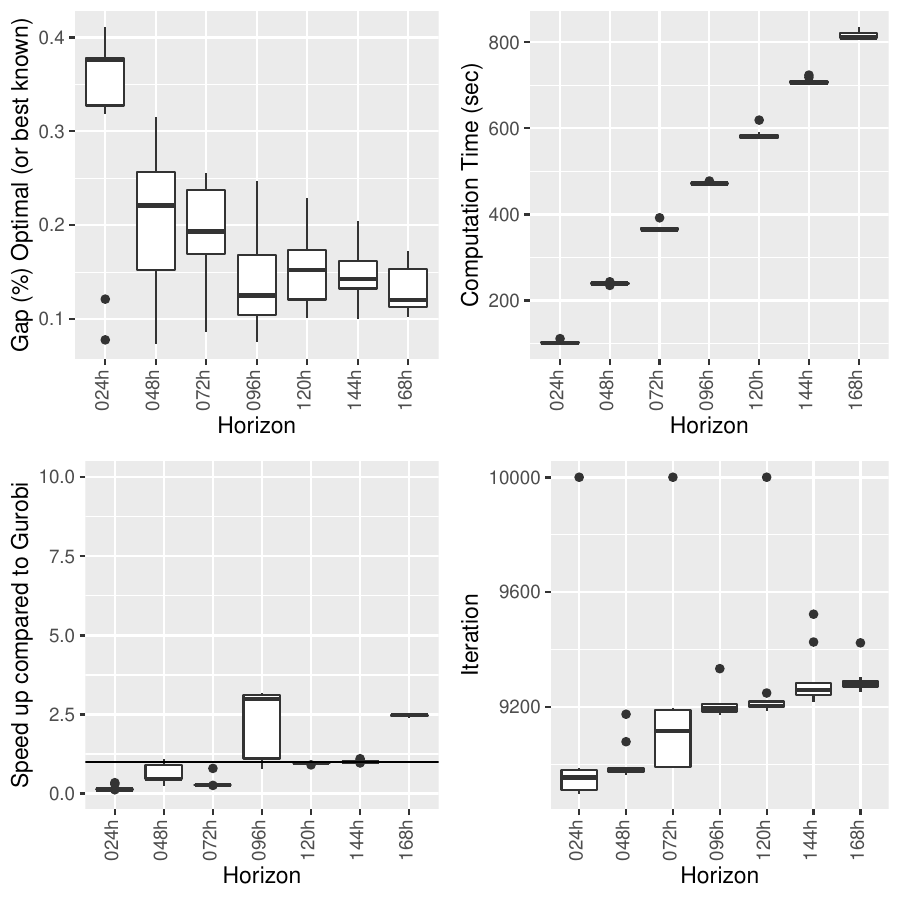}
    \caption{Additional experiments with the UC instance \textit{RTS54} with $\alpha=1.05$ and $m=5$}
    \label{fig:TransmissionResults}
\end{figure}

As our ADMM algorithm with the setting of $\alpha = 1.1$ and $m = 1$ cannot find a good solution for the instance RTS54, we did additional experiments with $\alpha$ set at 1.05 and m at 5 and presented these much-improved results in Figure~\ref{fig:TransmissionResults}. However, Gurobi seems just to be very successful at solving this particular instance. But if we consider the speedup, as the time horizon increases, we see that our algorithm gets better relative to Gurobi. This suggests that if we would increase the time horizon even more then it favors our algorithm.

\section{Discussion}

The results of our new ADMM algorithm for the UC looks promising in comparison with the competing state-of-the-art MIP approach. The advantages of the latter in general are that it can solve large programs, produce near optimal solutions, and can be easily extended. These three advantages also apply to our algorithm; the first two are demonstrated by our results and also our algorithm can easily be extended by implementing new subproblems as a QP or MIQP, as we have shown with the storage units.  Additionally, our method has the advantage that solving subproblems can also be improved by using more specific sophisticated algorithms. For example, our 1UC algorithm on average was almost 400 times faster than the MIQP approach \cite{wuijts2021improved}.  

Moreover, our algorithm seems to scale well with the number of time steps. The previous section showed that the time complexity in practice was almost linear. When solving instances with long time horizons, as is needed in large scale power system modeling, our method can find a high-quality solution in seconds while the MIP solver could not find any good solution even after 1 hour.

For future work this algorithm can be extended or improved in the following way:
\begin{itemize}
    \item Our algorithm currently does not include reserve requirements. Since reserve requirements are coupling constraints that includes decision variables from multiple components of the power system, we would need to relax these. Moreover, the addition of reserve requirements leads to additional Lagrange multipliers that we need to incorporate into our 1UC dynamic programming algorithm as was previously done by \cite{li2013dynamic}.
    \item Our algorithm currently does not include time dependent  start-up cost for the 1UC subproblems but this could be easily incorporated \cite{frangioni2006solving}.
    \item 	The subproblems can be modelled in a more complex but also more accurate way. For example, hydropower systems can be modelled in more detail  \cite{van2018large}, a more detailed transmission system (DC, AC) can be used, or a more robust transmission system that respect $n-1$ criteria can be modelled \cite{van2014dc}. Moreover, different power system assets such as demand response can easily be added as a subproblem.
\end{itemize}

\section{Conclusion}

In this paper, we presented a new ADMM algorithm that solves UC to near optimality. First, we transformed a standard UC formulation into an augmented Lagrangian formulation by relaxing the demand coupling constraints. Next, we solved the resulting problem with ADMM, an iterative algorithm. In each iteration we solved the subproblems induced by the relaxation of the demand coupling constraints. We used efficient solution methods for them. Especially, the use of our recently introduced fast algorithm for the 1UC subproblem enables many iterations in a limited time, so that ADMM can converge to a feasible solution of high quality. Morover, our algorithms is the first Langrangian beased method that finds feasible solutions for the case where ramping limits are included without requiring an ad-hoc repair heuristic. 

Computational experiments on a large set of UC instances showed that our algorithm produces high quality solutions. For instances with linear generation cost it outperformed the state-of-the-art MILP formulation, especially with longer time horizons. For the instances with quadratic cost our algorithms solved those instances significantly faster compared to the state-of-the-art MIQP formulation. As the time complexity of our ADMM grew almost linearly with the number of time steps, it could find high-quality solutions in seconds for instances with long time horizons, while the MIP solver could not find any good solution within one hour. 

Our ADMM, therefore, has potential to improve large scale power system modelling commonly based on UC with long time horizons, and planning of unit commitment in the short-term at high temporal resolution. Moreover, due to the decoupled nature of ADMM, the subproblem formulations can easily be adapted, for example to model elements of the power system in more detail or to improve the time complexity of solving the subproblems.

\section{Acknowledgement}
This work is part of the research programme ``Energie: Systeem 
Integratie en Big Data" with project number 647.003.005, which is financed by the Dutch Research Council (NWO).

\bibliography{sn-bibliography}

\appendix
\section{Variable Splitting and Gauss Seidel}

\subsection{Multi-block ADMM by Variable Splitting}\label{sec:MultiBlockSplitting}

Another way to decompose a problem with $n$ (sets) of variables is by variable splitting. First, we can rewrite the optimization problem of (\ref{eq:chapter3:basicObjective}) and (\ref{eq:chapter3:basicConstraint}) (or a problem with $n$ variables) to the following equivalent optimization problem that has additional variables $z_i$ that are copies of $p_i$:
\begin{align}
    \min \sum_{i=1}^n f_i(p_i)\\
    s.t. \nonumber\\
    D - \sum_{i=1}^n z_i = 0 \label{eq:chapter3:copiedDemandCoupling}\\
    p_i = z_i && \forall i \in \{1,\dots,n\} \label{eq:chapter3:copy}
\end{align}
 Afterwards, we can define the augmented Lagrangian by relaxing (\ref{eq:chapter3:copy}):
\begin{align}
    &    \mathcal{L}(\mathbf{p},\mathbf{z},\mathbf{\pi}) = \sum_{i=1}^n f_i(p_i) + \pi_i     (p_i - z_i) + \frac{\rho}{2} (p_i - z_i)^2  \\
    & \text{s.t. }     D - \sum_{i=1}^n z_i = 0 
\end{align}
The ADMM iterations are presented in Algorithm~\ref{alg:ADMMVariableSplitting}.
\begin{algorithm}[!htp]
\caption{Multi-block by Variable Splitting}
\label{alg:ADMMVariableSplitting}
\begin{algorithmic}[1]
\While{Stopping criteria have not been met}
\ForAll{$i = 1,\dots,n$}
 \State $ p_i^{k+1} \leftarrow \text{argmin}_{p_i} \{ f_i(p_i) + \pi_i^k p_i + \frac{\rho}{2}(p_i -  z_i^k )^2 \} $ 
 \EndFor
 \State $   \boldsymbol{z}^{k+1} \leftarrow \text{argmin}_{z_1,\dots, z_n } \{ \sum_{i=1}^n \pi_i^k (p_i^{k+1} -  z_i) + \frac{\rho}{2} (p_i^{k+1} - z_i)^2 \mid D  = \sum_{i=1}^n z_i \} $
\ForAll{$i = 1,\dots,n$}
 \State $ \pi_i^{k+1} \leftarrow  \pi_i^{k} + \rho  ( p_i^{k+1} - z_i^{k+1}) $
  \EndFor
 \State $k \leftarrow k +1$
\EndWhile
\end{algorithmic}
\end{algorithm}
The optimal values for any $p_i$ does not depend on the other production variables  because only the copied variables are coupled in (\ref{eq:chapter3:copiedDemandCoupling}). Therefore, at each iteration a similar subproblem is solved as when we used the Gauss-Seidel method. However, the difference is that the optimal value for all the copied variables must be found and more Lagrange multipliers must be updated. A main advantage of this method is that:
\begin{itemize}
    \item At its core it is still 2-block which has convergence guaranties.
    \item  The updates on $p_1$, $p_2$ ect. can be done in parallel instead of sequential because only the values from iteration $k$ and not from $k+1$ are used. 
\end{itemize}

There seems to be a disadvantage of introducing multiple additional variables and Lagrange multipliers and finding the optimal values of the copied variables at line $5$ in Algorithm~\ref{alg:ADMMVariableSplitting}. However, it is possible to achieve the same effect with only one copied variable and one Lagrange multiplier for each relaxed constraint.

We can rewrite the $z_i$ minimization step as follows:

\begin{align}
 \sum_{i=1}^n \pi_i^k (p_i^{k+1} - z_i) + \frac{\rho}{2} (p_i^{k+1} - z_i)^2 = \\
 \sum_{i=1}^n \frac{\rho}{2} \left[(p_i^{k+1} - z_i)^2 + 2 \frac{\pi_i^k}{\rho} \left( p_i^{k+1} - z_i\right) + \left(\frac{\pi_i^k}{\rho}\right)^2 - \left(\frac{\pi_i^k}{\rho}\right)^2 \right] = \\
  \sum_{i=1}^n \frac{\rho}{2} ((p_i^{k+1} - z_i) + \frac{\pi_i^k}{\rho})^2 - \frac{(\pi_i^k)^2}{2\rho} 
\end{align}

For $z_i$ to be optimal, we need to minimize 
subject to $\sum_i^n z_i=D$. We are minimizing a sum of quadratic differences $p_i^{k+1}  + \frac{\pi_i^k}{\rho} - z_i$. It is easy to show (by e.g the KKT conditions) that the minimum is attained if all the differences are equal, i.e., the difference would be the same for every $i$ and would be equal to the average difference, which implies
\begin{equation}
p_i^{k+1} + \frac{\pi_i^k}{\rho} - z_i = \frac{\sum_{j=1}^n (p_j^k + \frac{\pi_j^k}{\rho} - z_j)}{n} \label{eq:chapter3:avgDistance}
\end{equation}
\noindent
Since the total sum $\sum_i^n z_i$ must be equal to $D$ we can get $z_i^{k+1}$ by calculating:
\begin{align}
    z_i^{k+1} = \frac{D}{n} + p_i^{k+1} +  \frac{\pi_i^k}{\rho} - \frac{ \sum_{j=1}^n ({p_j^{k+1} +  \frac{\pi_j^k}{\rho}})}{n}
\end{align}
\noindent
which simplifies to:
\begin{align}
    z_i^{k+1} = \frac{D}{n} - \overline{p}^{k+1} + p_i^{k+1},  \label{eq:chapter3:ziAnalytical}
\end{align}
where $\overline{p}^{k+1} = \sum_{i=1}^n p_i^{k+1}$

If all  multipliers $\pi_i$ start with the same value at the beginning of the iteration then the update (Algorithm~3 line 7) would be the same for every multiplier due to (\ref{eq:chapter3:avgDistance}). This implies that all multipliers have the same value and can be replaced by a single multiplier. Moreover, since the optimal value of $z_i^{k+1}$ has an analytical form (\ref{eq:chapter3:ziAnalytical}), we can simplify Algorithm~3 to Algorithm~4.

\begin{algorithm}[!htp]
\caption{Optimal Exchange}
\label{alg:exchange}
\begin{algorithmic}[1]
\While{Stopping criteria have not been met}
\ForAll{$i = 1,\dots,n$}
 \State $ p_i^{k+1} \leftarrow \text{argmin}_{p_i} \{ f_i(p_i) + \lambda^k p_i + \frac{\rho}{2}(p_i -  \frac{D }{n} +\overline{p}^{k}   - p_i^k  )^2 \} $ 
 \EndFor

\ForAll{$i = 1,\dots,n$}
 \State $\lambda^{k+1} \leftarrow \lambda^k +  \rho (\overline{p}^{k+1}-\frac{D }{n})$
  \EndFor
 \State $k \leftarrow k +1$
\EndWhile
\end{algorithmic}
\end{algorithm}

In fact Algorithm~4 is called the \textit{optimal exchange} problem which in turn is a special case of the \textit{sharing problem} \cite{boyd2011distributed}. \textit{Optimal exchange} is similar to  Gauss-Seidel, they both have a quadratic constant that is based on the residual load. However, they differ with respect to convergence, as \textit{optimal exchange} is proven to be convergent, but Gauss Seidel not necessarily. Moreover, finding the optimal value of $p_1,\dots,p_n$ can be done in parallel with \textit{optimal exchange} but not with Gauss Seidel. 

\subsection{Gauss Seidel Compared to Variable Splitting}

\begin{figure}[!htbp]
    \centering
    \includegraphics[scale=0.9]{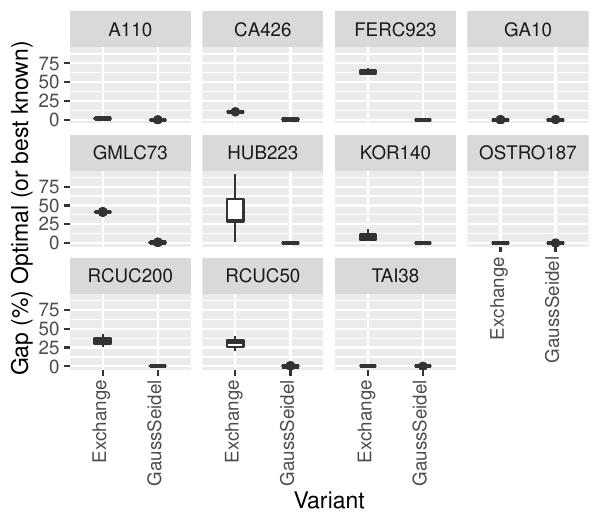}
    \caption{Gap when using the Exchange or Gauss Seidel ADMM variation. Time horizon for experiments were $24$ hours. All instances were repeated 10 times.}
    \label{fig:exchangeResults}
\end{figure}

An unexpected result is that the Gauss Seidel outperformed the Exchange ADMM variation by a large margin. Gauss Seidel converged faster and found heuristic solutions of better quality.  For most instances the solution found by the Exchange algorithm was far from optimal (Figure~\ref{fig:exchangeResults}), while Gauss Seidel found high-quality solutions. Therefore, for the remainder of the results and figures we only show the results of our Gauss Seidel variant.

\end{document}